\documentclass[a4paper,12pt]{amsart}

\usepackage{amsmath,latexsym,amscd}
\usepackage{amssymb} 
\usepackage[mathscr]{eucal} 
\usepackage[all]{xy}

\newfont{\sheaf}{eusm10 scaled\magstep1}

\title{The Noether inequality for smooth minimal 3-folds}
\author{Fabrizio Catanese, Meng Chen and De-Qi Zhang}
\address{\rm Lehrstuhl Mathematik VIII,
Universit$\ddot{\text{a}}$t Bayreuth, D-95440, BAYREUTH, Germany}
\email{Fabrizio.Catanese@uni-bayreuth.de}
\address{\rm School of Mathematical Sciences, Fudan University,
Shanghai, 200433, PR China} \email{mchen@fudan.edu.cn}
\address{\rm Department of Mathematics, National University of Singapore,
2 Science Drive 2, Singapore 117543, Singapore}
\email{matzdq@nus.edu.sg}
\thanks{
The present research took place in the
framework of the DFG Schwerpunkt "Globale Methode in der komplexen Geometrie", and of
the  joint Chinese-German project "Komplexe Geometrie" supported by the DFG and the NSFC.
Chen was supported by the National Natural Science Foundation of China.
Zhang was supported by an Academic Research Fund of NUS}

\newcommand{\roundup}[1]{\ulcorner{#1}\urcorner}

\newcommand{\Z}{\ensuremath{\mathbb{Z}}}

\newcommand{\hol}{\ensuremath{\mathcal{O}}}

\newcommand{\PP}{\ensuremath{\mathbb{P}}}

\newcommand{\ra}{\ensuremath{\rightarrow}}

\newcommand{\om}{\omega}

\newcommand{\Ga}{\Gamma}
\newcommand{\De}{\Delta}

\newtheorem{thm}{Theorem}[section]
\newtheorem{lem}[thm]{Lemma}
\newtheorem{cor}[thm]{Corollary}

\newtheorem{claim}[thm]{Claim}
\theoremstyle{definition}

\newtheorem{setup}[thm]{}
\newtheorem{question}[thm]{Question}

\newtheorem{conj}[thm]{Conjecture}
\newtheorem{rem}[thm]{Remark}
\theoremstyle{remark}

\def\eea{\end{eqnarray*}}
\def\bea{\begin{eqnarray*}}
\def\hol{{\mathcal{O}}}

\begin{document}
\begin{abstract} Let $X$ be a smooth projective minimal 3-fold of
general type. We prove the sharp inequality
$$K_X^3\ge \frac{2}{3}(2p_g(X)-5),$$
an analogue of the classical Noether inequality for
algebraic surfaces of general type.
\end{abstract}
\maketitle
\pagestyle{myheadings} \markboth{\hfill F. Catanese, M. Chen and
D-Q. Zhang\hfill}{\hfill The Noether inequality \hfill}
\section{\bf Introduction}

In the 1980's M. Reid,  observing the
importance of the Noether inequality: $K^2\ge 2p_g-4$ for
surfaces of general type, asked the following question

\begin{question}\label{Q}(M. Reid) What is the 3-dimensional version of
Noether's  inequality?
\end{question}

Question \ref{Q} is obviously a very important aspect of threefold
geography, just like the well known Miyaoka-Yau inequality. 
  There have
been already several works dedicated to the above question:
\begin{quote}
\item$\bullet$ M. Kobayashi (1992, \cite{kobayashi1992}) studied
Gorenstein minimal 3-folds of general type and found an infinite
number of examples (Proposition 3.2 in \cite{kobayashi1992})
satisfying the equality:

\noindent (1.1)\hskip2cm  $K^3=\frac{4}{3}p_g-\frac{10}{3}.$

\item$\bullet$ M. Chen (2004, \cite{chenjmsj2004}) gave effective
Noether type  inequalities for arbitrary minimal 3-folds of
general type.

\item$\bullet$ M. Chen (2004, \cite{chenmrl}) answered Question
\ref{Q} under the assumption that the 3-fold $X$ is smooth with an
ample canonical line bundle, proving the sharp inequality:
$K^3\ge \frac{4}{3}p_g-\frac{10}{3}$.

[{\tiny In the above three items, $K^3:=K_X^3$ is the canonical
volume and $p_g:=p_g(X)$ is the geometric genus of $X$.}]

\end{quote}

In this paper, we will generalize the main theorem of
\cite{chenmrl}. The aim is to answer Question \ref{Q} under a
weaker condition:
\begin{thm}\label{main}
Let $X$ be a smooth projective minimal 3-fold of general type.
Then the sharp Noether inequality:
$$K_X^3\ge \frac{2}{3}(2p_g(X)-5)$$
holds.
\end{thm}

\begin{rem} The inequality in Theorem \ref{main} is sharp because
of M. Kobayashi's interesting examples (cf. Equation (1.1)).
\end{rem}

As an application of our results, we present the following
corollary which gives a classification of 3-folds of general type
with small "slope" $K^3/p_g$:

\begin{cor}\label{application} Let $X$ be a projective minimal
(i.e., $K_X$ is nef) Gorenstein 3-fold of general type with canonical
singularities. Assume $K_X^3<\frac{7}{5}p_g(X)-2$. Then $X$ is
canonically fibred by curves of genus 2.
\end{cor}

The assumption in Corollary \ref{application} is not empty again
because of M. Kobayashi's examples.

\begin{setup}{\bf The set up.}\label{notation}
Let $X$ be a projective minimal  Gorenstein 3-fold of general type
with canonical singularities. According to the work of M. Reid
\cite{reid1983} and Y. Kawamata (Lemma 5.1 of)
\cite{kawamata1988}, there is a minimal model $Y$ with a
birational morphism $\nu:Y\longrightarrow X$ such that
$K_Y=\nu^*(K_X)$ and such that $Y$ is factorial with at worst terminal
singularities. Thus we may always assume that $X$ is factorial
with only (necessarily finitely many)
  terminal singularities. Observing that $K_X^3\ge 2$ (see
\ref{even} below), the inequality in Theorem \ref{main} is
automatically true whenever $p_g(X)\le 4$. So the essential argumentation
takes place when $p_g(X)$ is bigger and we are led to study the
canonical map $\Phi := \Phi_{|K_X|}$ as in the two dimensional case.

Take a birational modification $\pi: X'\longrightarrow X$, which
exists by Hironaka's big theorem, such that:

(1) $X'$ is smooth;

(2) the movable part of $|K_{X'}|$ is base point free;

(3) $\pi^*(K_X)$ is supported by a normal crossings divisor (so that
we are in a position to apply the Kawamata-Viehweg vanishing theorem
\cite{kawamata1982,viehweg1982}).

We will fix some notation below. Denote by $g$ the composition
$\Phi\circ\pi$. So $g: X'\longrightarrow W'\subseteq{\mathbb P}^N$
is a morphism. Let $g: X'\overset{f}\longrightarrow B
\overset{s}\longrightarrow W'$ be the Stein factorization of $g$
(thus $B$ is normal and $f$ has connected fibers). We can
write:
$$K_{X'}=\pi^*(K_X)+E_{\pi}=M+Z',$$
where $M$ is the movable part of $|K_{X'}|$, $Z'$ the fixed part
and $E_{\pi}$ an effective divisor which is a linear combination
of distinct exceptional divisors. We may also write:
$$\pi^*(K_X)=M+E',$$
where $E'=Z'-E_{\pi}$ is an effective divisor. On $X$, one may
write $K_X\sim N+Z$ where $N$ is the movable part and $Z$ the
fixed part. So
$$\pi^*(N)=M+\sum_{i=1}^sd_iE_i$$
with $d_i>0$ for all $i$. The above sum runs over all those
exceptional divisors of $\pi$ that lie over the base locus of $M$.
On the other hand, one may write $E_{\pi}=\sum_{j=1}^te_jE_j$
where the sum runs over {\it all} exceptional divisors of $\pi$.
One has $e_j>0$ for all $1 \le j \le t$ because $X$ is terminal.
Apparently, one has $t\ge s$.

Set $d:=\dim(B)$. We say that $X$ {\it is canonically fibred by
surfaces} if $d=1$. Under this situation, we have an induced
fibration $f:X'\longrightarrow B$ onto a smooth curve $B$. Denote
by $b:=g(B)$ the geometric genus of $B$.
\end{setup}

{\bf Notations}
\begin{tabbing}
  \= aaaaaaaaaaaaaaaaa \= bbbbbbbbbbbbbbbbbbbbbbbbbbbbbbb \kill

\> $K^3$    \> the canonical volume of a 3-fold in question\\
\> $p_g = h^0 (\hol(K))$    \> the geometric genus\\
\> $q(V) =h^1 (\hol_V)$    \> the irregularity of $V$\\
\> $h^2({\mathcal O}_V)$ \> the second irregularity of a 3-fold
$V$\\
\> $\chi({\mathcal O_V})$ \> the Euler Poincare characteristic of $V$\\
\> $(K^2, p_g)$ \> invariants of a minimal surface of general
type\\
\> $g(B)$ \> the genus of a curve $B$\\
\>$\equiv$     \> numerical equivalence\\
\> $\sim$      \> linear equivalence\\
\> $\roundup{\cdot}$  \> the round up of $\cdot$
 ($\roundup{x}:= min \{ n\in \Z| n
\geq x \}$)\\
\> $D|_S$ \> the restriction of the divisor $D$ to $S$\\
\> $D\cdot C$ \> the intersection number of a divisor $D$ with a curve $C$\\
\end{tabbing}

\section{\bf Reduction to the surface case and the lower bound of $K^3$}
\begin{setup}\label{even}{\bf $K^3$ is even.}
Suppose that $D$ is any divisor on a smooth 3-fold $V$. The
Riemann-Roch theorem (cf. appendix in Hartshorne's book
\cite{Hart}) gives:
$$\chi({\mathcal O}_{V}(D))=\frac{D^3}{6}-\frac{K_V\cdot D^2}{4}+
\frac{D\cdot(K_V^2+c_2(V))}{12}+\chi({\mathcal O}_{V}).$$
A direct calculation shows that
$$\chi({\mathcal O}_{V}(D))+\chi({\mathcal O}_{V}(-D))=\frac{-K_V\cdot
D^2}{2}+2\chi({\mathcal O}_{V})\in{\mathbb Z}.$$ Therefore,
$K_V\cdot D^2$ is an even number.

Now let $X$ be a projective minimal Gorenstein 3-fold of general
type. Denote by $\nu:V\longrightarrow X$ a smooth birational
modification. Let $D$ be any divisor on X. Then $K_X\cdot
D^2=K_V\cdot (\nu^*D)^2$ is even. Especially, $K_X^3$ is even and
positive.
\end{setup}

\begin{setup}\label{known}{\bf Known results.} Let $X$ be a projective
minimal factorial 3-fold of general type with terminal
singularities. The following Noether type 
inequalities have already been established, where $d = \dim \Phi_{|K_X|}(X)$.

\begin{quote}
\item $\bullet$ if $d=3$, then $K_X^3\ge 2p_g(X)-6$ 
 (cf.  M. Kobayashi's Main Theorem in \cite{kobayashi1992});

\item $\bullet$ if $d=2$, then $K_X^3\ge
\roundup{\frac{2}{3}(g-1)}(p_g(X)-2)$  (cf.  Chen's Theorem 4.1(ii)
in \cite{chenjmsj2004}), where $g$ is the genus of a general fiber
of the induced fibration $f:X'\longrightarrow B$; if furthermore
$X$ is smooth, then $K_X^3\ge \frac{2}{3}(2p_g(X)-5)$  (cf.  Chen's
Theorem 4.3 in \cite{chenjmsj2004});

\item $\bullet$ if $d=1$ and the general fiber $S$ of the induced
fibration $f:X'\longrightarrow B$ is not a surface of type
$(K^2, p_g)=(1,2)$, then $K_X^3\ge 2p_g(X)-4$   (cf. Chen's
Theorem 4.1(iii) in \cite{chenjmsj2004}).
\end{quote}
\medskip

In order to prove Theorem \ref{main}, we have to treat the remaining
case (in the above third item) where $X$ is canonically fibred by
surfaces of type $(1,2)$. Note that Theorem \ref{main} was proved
in \cite{chenmrl} only under the stronger assumption of $K_X$ being
ample. Assuming only  the nefness of $K_X$, we can see that the method in
\cite{chenmrl} is no longer effective and the situation could be
more complicated. It is the aim of this paper to overcome
this obstacle and prove our Theorem \ref{main}.
\end{setup}

The rest of this section is devoted to deducing several key
inequalities through the ${\mathbb Q}$-divisor method.

\begin{setup}\label{key}{\bf Key inequalities.}
Keep the same notation as in \ref{notation} and assume that $K_X$
is nef and big. Suppose, from now on, $d=1$ and $p_g(X)\ge 3$. We
have an induced fibration $f:X'\longrightarrow B$. Denote by $S$ a
general fiber of $f$. Let $\sigma:S\longrightarrow S_0$ be the
contraction onto the minimal model. Suppose $(K_{S_0}^2,
p_g(S_0))=(1,2)$.

By Lemma 4.5 of \cite{chenjmsj2004}, we have two cases exactly:
$$q(X)=b=1\ \ \text{and}\ \ h^2({\mathcal O}_X)=0,$$
$$q(X)=b=0\ \ \text{and}\ \ h^2({\mathcal O}_X)\le 1.$$

One may write $M=\sum_{i=1}^aS_i$ as a disjoint union of distinct
smooth fibers of $f$, where $a=p_g(X)-1$ if $b=0$, or
$a=p_g(X)$ otherwise. Noting that $\pi^*(K_X)_{|S} \le K_S$ is a nef and
big Cartier divisor and that $\sigma^*(K_{S_0})$ is the positive
part of the Zariski decomposition of $K_S$, so
$\pi^*(K_X)_{|S}^2=\sigma^*(K_{S_0})^2=1$, and
$\pi^*(K_X)_{|S}\sim \sigma^*(K_{S_0})$ by the uniqueness of the
Zariski decomposition. According to the construction of $\pi$, we
know that $E'_{|S}\sim \pi^*(K_X)_{|S}$ is a normal crossing
divisor for a general fiber $S$.

Now let us assume $\alpha_3\in (0,1)$ be a real number such that
$$h^0(S, K_S+\roundup{\alpha E'_{|S}})\ge 3$$
for all $\alpha>\alpha_3$. We may now write $a$ as $a=m_2+m_3+1$, where 
$m_2, m_3 $ are non-negative integers and
$$\frac{a-m_3}{a}>\alpha_3.$$
Such integers exist: for instance, one
may take $m_3=0$ and $m_2=a-1$. What we will show in next sections
is that we can find a nontrivial decomposition of $a$, i.e., with $ m_3 > 0$.

Once we have the above setting, we may deduce an interesting
inequality as follows. Write
$$M\sim
S_0+\sum_{i=1}^{m_2}S_{2,i}+\sum_{j=1}^{m_3}S_{3,j}.$$
Since
$$\pi^*(K_X)-\sum_{j=1}^{m_3}S_{3,j}-\frac{m_3}{a}E'\equiv
(1-\frac{m_3}{a})\pi^*(K_X)$$ is nef and big and has normal crossings,
the Kawamata-Viehweg vanishing theorem (\cite{kawamata1982,
viehweg1982}) yields
$$H^1(X',
K_{X'}+\roundup{\pi^*(K_X)-\sum_{j=1}^{m_3}S_{3,j}-\frac{m_3}{a}E'})=0$$
and hence the exact sequence:
\begin{align*}
0&\longrightarrow
H^0(X',K_{X'}+\roundup{\pi^*(K_X)-\sum_{j=1}^{m_3}S_{3,j}
-\frac{m_3}{a}E'})\\
  &\longrightarrow H^0(X',
K_{X'}+\roundup{\pi^*(K_X)-\frac{m_3}{a}E'})\\
&\longrightarrow\oplus_{j=1}^{m_3}
H^0(S_{3,j},K_{S_{3,j}}+\roundup{(1-\frac{m_3}{a})E'}_{|{S_{3,j}}})
\longrightarrow
0.
\end{align*}
In the above sequence, we obviously have
$$\roundup{(1-\frac{m_3}{a})E'}_{|{S_{3,j}}}\ge
\roundup{(1-\frac{m_3}{a})E'_{|{S_{3,j}}}}$$ and
$$(1-\frac{m_3}{a})E'_{|{S_{3,j}}}\equiv
\frac{a-m_3}{a}\pi^*(K_X)_{|S_{3,j}}.$$ So one has
$$h^0(S_{3,j}, K_{S_{3,j}}+\roundup{(1-\frac{m_3}{a})E'_{|{S_{3,j}}}})
\geq 3$$
for sufficiently general $S_{3,j}$ as a fiber of $f$ by our
definition of $\alpha_3$. The above sequence then gives the
inequality
\begin{align*}
\tag{2.1}&h^0(X', K_{X'}+\roundup{\pi^*(K_X)-\frac{m_3}{a}E'})\\
&\ge
h^0(K_{X'}+\roundup{\pi^*(K_X)-\sum_{j=1}^{m_3}S_{3,j}-\frac{m_3}{a}E'})
+3m_3.
\end{align*}
It is obvious that one has
\begin{align*}
&h^0(K_{X'}+\roundup{\pi^*(K_X)-\sum_{j=1}^{m_3}S_{3,j}-\frac{m_3}{a}E'})\\
&\ge
h^0(K_{X'}+\roundup{\pi^*(K_X)-\sum_{j=1}^{m_3}S_{3,j}
-\frac{m_2+m_3}{a}E'}).
\end{align*}
Similarly, because
$$\pi^*(K_X)-\sum_{i=1}^{m_2}S_{2,i}-\sum_{j=1}^{m_3}S_{3,j}
-\frac{m_2+m_3}{a}
E'\equiv
\frac{1}{a}\pi^*(K_X)$$ is nef and big and with normal crossings, the
vanishing theorem gives
$$H^1(K_{X'}+\roundup{\pi^*(K_X)-\sum_{i=1}^{m_2}S_{2,i}
-\sum_{j=1}^{m_3}S_{3,j}
-\frac{m_2+m_3}{a}E'})=0.$$ So we have the following exact
sequence:
\begin{align*}
0&\longrightarrow H^0(X',
K_{X'}+\roundup{\pi^*(K_X)-\sum_{i=1}^{m_2}S_{2,i}-
\sum_{j=1}^{m_3}S_{3,j}-
\frac{m_2+m_3}{a}E'})\\
&\longrightarrow
H^0(X',K_{X'}+\roundup{\pi^*(K_X)-\sum_{j=1}^{m_3}S_{3,j}
-\frac{m_2+m_3}{a}E'})\longrightarrow\\
&\oplus_{i=1}^{m_2}H^0(S_{2,i},
K_{S_{2,i}}+\roundup{\frac{a-m_2-m_3}{a}E'}|_{S_{2,i}})\longrightarrow
0.
\end{align*}
The above exact sequence gives
$$h^0(S_{2,i}, K_{S_{2,i}}+\roundup{\frac{a-m_2-m_3}{a}E'}|_{S_{2,i}})
\ge p_g(S_{2,i})=2$$
and
\begin{align*}
\tag{2.2}&h^0(X',K_{X'}+\roundup{\pi^*(K_X)-\sum_{j=1}^{m_3}S_{3,j}
-\frac{m_2+m_3}{a}E'} )\\
&\ge
h^0(K_{X'}+\roundup{\pi^*(K_X)-\sum_{i=1}^{m_2}S_{2,i}
-\sum_{j=1}^{m_3}S_{3,j}-\frac{m_2+m_3}{a}E'})+2m_2.
\end{align*}
We shall go on studying the group $$H^0(X',
K_{X'}+\roundup{\pi^*(K_X)-\sum_{i=1}^{m_2}S_{2,i}
-\sum_{j=1}^{m_3}S_{3,j}-
\frac{m_2+m_3}{a}E'}).$$
Apparently, it is slightly bigger than
$H^0(X', K_{X'}+S_0)$.

We set $\delta:=2-h^2({\mathcal O}_X).$ By looking at the exact
sequence:
$$0\longrightarrow {\mathcal
O}_{X'}(K_{X'})\longrightarrow{\mathcal
O}_{X'}(K_{X'}+S_0)\longrightarrow {\mathcal
O}_{S_0}(K_{S_0})\longrightarrow 0,$$ one has
\begin{align*}\tag{2.3} h^0(K_{X'}+S_0)\ge p_g(X)+\delta.
\end{align*}
Combining the above inequalities (2.1)$\sim$(2.3), we have
\begin{align*}\tag{2.4}
P_2(X)=h^0(K_{X'}+\pi^*(K_X))\ge 3m_3+2m_2+p_g(X)+\delta.
\end{align*}
Applying Reid's plurigenus formula (see the last section of
\cite{YPG} and Lemma 8.3 of \cite{Miyaoka}):
$$P_2(X)=\frac{1}{2}K_X^3-3\chi({\mathcal
O}_X)=\frac{1}{2}K_X^3-3(1-b+h^2({\mathcal O}_X)-p_g(X)),$$ we get
the Noether type inequality:
\begin{align*}\tag{2.5}
K_X^3\ge 6m_3+4m_2-4p_g(X)+4h^2({\mathcal O}_X)-6b+10.
\end{align*}
\end{setup}

\begin{setup}\label{problem}{\bf A problem on surfaces.}
As we have seen, the general fiber $S$ has the invariants
$(K_{S_0}^2, p_g(S_0))=(1,2)$. We have a divisor $E'_{|S}\sim
\sigma^*(K_{S_0})$ which has normal crossings. So there is a divisor
$D_0\in |K_{S_0}|$ with $E'_{|S}=\sigma^*(D_0)$. We expect to find a
real number $\alpha_3\in (0,1)$ such that $h^0(S,
K_S+\roundup{\alpha E'_{|S}})\ge 3$ for all $\alpha>\alpha_3$.
Furthermore we hope $\alpha_3$ to be as small as possible.
\end{setup}

\section{\bf The rounding up problem for (1,2) surfaces}

Assume that $Y$ is the canonical model of a surface of general
type with $p_g(Y)=2, K^2_Y=1$,  that $\tau : S_0 \ra Y$ is its
minimal model, and finally that $f : S \ra S_0$ is a sequence of
point blow ups.

We set up the following notation and assumptions:
\begin{itemize}
\item $\Ga \subset Y $ is a canonical divisor \item $D$ is the
full transform $ \tau ^* ( \Ga )$ \item we assume that $f^* (D)$
is a normal crossing divisor \item for $ t \in (0,1)$ we consider
the round up divisor $ \De _t : =\roundup { t f^* (D)}$
\end{itemize}

\begin{rem}
Observing that since $ H^1 (\hol_Y) =  H^1 (\hol_S) =H^1 (K_S) =
0$ (cf. \cite{BPV}), one has $ h^0 ( K_S + \De _t) = p_g(S) + h^0
( \om_{ \De _t})$ where $\om_{ \De _t}:={\mathcal O}_{\De _t}(K_S
+ \De _t)$.
\end{rem}

\begin{thm}\label{cat} Assume that $p_g(Y)=2, K^2_Y=1$, that $f^* (D)$ is a
normal crossings divisor, and that $ 3/10 < t$. Then $ h^0 ( K_S +
\De _t) = 2 + h^0 ( \om_{ \De _t}) \geq 3.$
\end{thm}

\begin{proof}

1) Since $ K_Y^2 = 1$, and $ K_Y$ is ample, $ \Gamma$ is
irreducible. (Note that $|K_Y|$ has one smooth and simple base point and the general
member of $|K_Y|$ is a smooth curve of genus 2 (cf. page 225 in
\cite{BPV}). It is well known and easy to show that $Y$ is a hypersurface of degree $10$
in the weighted projective space $\PP(1,1,2,5)$, so $Y$ is a  finite double cover 
of $\PP(1,1,2)$ and  the involution
$\sigma$ on
$Y$ induced by the hyperelliptic involution of those genus 2 curves has exactly one isolated
fixed point -- the base point of $|K_Y|$.
We shall also denote by the same symbol $\sigma$ its lift to a biregular involution
on $S_0$, observing that again there is exactly one isolated
fixed point -- the base point of $|K_{S_0}|$. 

The quotient
$Q_2 = Y/\langle \sigma \rangle  = \PP(1, 1, 2)$
is  isomorphic to a quadric cone in
${\mathbb P}^3$ and $\Gamma$ is isomorphic to a double cover of
$\PP^1$ branched in a point $ P_{\infty}$ and in a disjoint
sub-scheme of length $5$ (cf. \cite{BPV}, page 231,
a construction due to Horikawa).

2)  Observe that if $ D \geq D'$,  and $ \De '_t : =\roundup { t
f^* (D')}$, then $ h^0 ( K_S + \De _t) \geq  h^0 ( K_S + \De'
_t)$.

3) Set $ K : = K_{S_0}$. Write $ D = \tilde {\Ga } + \tilde{Z}$,
where $ \tilde {\Ga } $ is the strict transform of $\Ga$.
Thus  $ \tilde {\Ga } \cdot K =
1, \tilde{Z} \cdot K = 0$. Since $\Ga$ is
a Cartier divisor, it follows that the support of $ \tilde{Z}$ is
a union of the support of certain fundamental cycles $Z_i$
corresponding to the rational double points $ P_i \in X$ such that
$ P_i \in \Ga$, and moreover $ \tilde{Z} = \sum_i \tilde{Z}_i$,
where $ \tilde{Z}_i \geq Z_i$.

4) If we take an effective decomposition $ D = D' + W$, where $ D'
\cdot K =1$, then $ (D')^2 = D' \cdot ( K - W) = 1 - D' \cdot W
\leq -1 $, since a canonical curve is 2-connected (\cite{BPV}, VII (6.2)).

5) If $Z' \cdot K = 0$, and $Z'$ is (effective and) reduced, then
$ (Z')^2 = - 2 k$, where $k$ is the number of connected components
of $Z'$. In fact, it suffices to prove the formula for $Z'$
connected, but $Z'$ is contained in a fundamental cycle, and
corresponds therefore to a rational subtree of the Dynkin diagram.
Thus, if $n$ is the number of edges of the subtree,
then $ (Z')^2 = - 2 (n+1) + 2 n = -2.$

We pass now to the strategy of proof:

{\bf [S1]} if the arithmetic genus $p (\tilde {\Ga } ) \geq 1$ then we
pick $ D' = \tilde {\Ga } $ (see point 2)).

Observe now that $p (\tilde {\Ga } ) \geq 1$ is equivalent, since
$ \tilde {\Ga } \cdot K = 1$, to  $ \tilde {\Ga } ^2 \ge -1$, or to
$ D =   \tilde {\Ga } $, in view of  4). If the first strategy is
not allowed, this means that $ \tilde {\Ga } ^2 = -3$, and $
\tilde {\Ga } \cong \PP^1.$

     If $ \tilde {\Ga } \cong \PP^1$ we consider the reduced divisor $
\tilde {\Ga } + Z'_i$, where $Z'_i = (Z_i)_{red}$ is the reduced
curve corresponding to one of the divisors $\tilde{Z}_i$ appearing
in 3). By 5) and 4) it follows that the odd number $( \tilde {\Ga
} + Z'_i)^2 = -5 + 2 ( \tilde {\Ga } \cdot  Z'_i)$ equals $ -3$ or
$-1$, accordingly $ ( \tilde {\Ga } \cdot  Z'_i) = 1 \ {\rm or} \
2.$

{\bf [S2]} If $ ( \tilde {\Ga } \cdot  Z'_i) = 2$, there are four cases:

{\bf [S2.1]} $  \tilde {\Ga } +   Z'_i $ is a normal crossing divisor
(of arithmetic genus $1$), and we pick $  D' = \tilde {\Ga } +
Z'_i $.

{\bf [S2.2]} $  \tilde {\Ga }$ is tangent to a smooth (-2)-curve $ A
\subset   Z'_i $, and then we take $ D' = \tilde {\Ga } + A $.

{\bf [S2.3]} A fundamental cycle $Z_1 < \tilde{Z}$ is of type $A_4$
and $\tilde{\Ga}$ passes through the central point transversally.
Take $D' = D$ (see the claim below).

{\bf [S2.4]} A fundamental cycle $Z_1 < \tilde{Z}$ is of type $A_2$
and $\tilde{\Ga}$ passes through the central point transversally.
Take $D' = \tilde{\Ga} + Z_2$ (see the claim below).

\begin{claim}
(1) Cases [S2.1] -- [S2.4] are the only possible cases if [S2] holds.

\par
(2) In Case [S2.3], one has $K_{S_0} \sim D = \tilde{\Ga} + \tilde{Z}$ with
$\tilde{Z} = A_1 + 2A + 2A' + A_4$, so that $Z_1 =  A_1 + A + A' + A_4$
is a fundamental cycle of type $A_4$.

\par
(3) In Case[S2.4], there is another fundamental cycle $Z_2 < \tilde{Z}$ 
of type $A_m$ which together with $\tilde{\Ga}$ forms a rational
loop (of arithmetic genus 1).

\end{claim}
\begin{proof} {\it (of the claim)} If $  \tilde {\Ga } +   Z'_i $ is not a
normal crossings divisor, then, the intersection number being $2$,
either [S2.2] holds or  $  \tilde {\Ga } $ meets $ Z'_i $ at a
singular point $P$ where two components $ A, A'$ meet, and all
intersections are transversal. We observed that on $ S_0$ we have a
canonical  biregular involution $\sigma$, induced from the
hyperelliptic involution on the (genus two) canonical curves.

$P$ is then a fixed point for the involution $\sigma$, which has only the
point lying over $P_{\infty}$ as isolated fixed point. Since $P$
lies in a fundamental cycle, $P$ is a different point than the
above isolated fixed point. So there is a $\sigma$-fixed curve $C$ (on $S_0$)
through $P$.
If both $A, A'$ are $\sigma$-stable, then the action $\sigma_*$
on the tangent space at $P$ will have three eigenvectors
(along $A, A', \tilde{\Ga}$) and hence it equals $(-1) id$,
contradicting the fact that $P$ is not an isolated $\sigma$-fixed point.

Thus $\sigma$ must interchange $A$ and $A'$.

Let $\tilde{Z}_1$ contain $A, A'$. Then $\sigma$ acts on the graph
of $\tilde{Z}_1$ fixing $P = A \cap A'$. So $\tilde{Z}_1$ is of Dynkin
type $A_{2n}$ ($n \ge 1$) and $P$ is the central point of $\tilde{Z}_1$.
Therefore, $A, A'$ are the inverse images in the double cover
$S_0 \rightarrow Q_2$ of the last exceptional curve of the blow up of
a singular
point $P'$ of the branch curve $B$ on $Q_2$.
Indeed, $P' \in B$ is a cusp of type $(2, 2n+1)$ with fibre $F$
the only tangent at $P' \in B$. By point 1)  follows that
$5 \ge (F . B)_{P'} = 2n+1$.
Thus $n = 1, 2$. This proves the first assertion.

The second assertion follows from  point 1).

Concerning assertion (3), by  point 1) and observing that
$\tilde{\Ga} \cong \PP^1$, our
$F$ has one further intersection
point $P_2$ with $B$, with $(F . B)_{P_2} = 2$, and with
$P_2$ a singular point for $B$ of type $A_n$. Then  assertion (3) follows.
\end{proof}

6) If strategies [S1] and [S2] are both not allowed, this means
that $ \tilde {\Ga } \cong \PP^1$, and that $ ( \tilde {\Ga }
\cdot Z'_i) = 1$ for each $i$.

7) Consider the intersection number $\tilde {\Ga } \cdot \tilde
{Z_i}$  which equals $ K \cdot  \tilde {Z_i} - ( \tilde {Z_i})^2 =
- ( \tilde {Z_i})^2 $ as is therefore a strictly positive even
number. By 4), since $(\tilde {\Ga } + \tilde {Z_i})^2 = (\tilde
{\Ga } )^2 - (\tilde {Z_i})^2 \le -1$ this number equals  $2$,
or $\tilde {Z_i} = \tilde {Z}$, in which case we get $4$
(indeed, note that $1 = \tilde{\Ga} \cdot (\tilde{\Ga} + \tilde{Z})$).

8) Assume still that strategies [S1] and [S2] are both not
allowed, thus $ 1 = K^2 = \tilde {\Ga }^2 + \sum_i  (2(\tilde {\Ga
} \cdot \tilde {Z_i}) +  \tilde {Z_i}^2) $, which, by 7), equals $
-3 + \sum_i  2$ if there is more than one fundamental cycle.
Therefore we conclude that $\tilde {\Ga }$ intersects precisely
one or two fundamental cycles, and in the former case $\tilde {\Ga
} \cdot \tilde {Z_i} = 4.$

9) Let us consider first the case where there are two fundamental
cycles intersecting $ \tilde {\Ga }$, and observe the following
inequalities: $2= \tilde {\Ga } \cdot \tilde {Z_i} \geq \tilde
{\Ga } \cdot Z_i
     \geq \tilde {\Ga } \cdot Z'_i = 1$, and write $ \tilde {Z_i} =
Z_i + W_i$.  We have $ (\tilde {\Ga } + \tilde {Z_i}) \cdot Z_i =
0 = \tilde {\Ga } \cdot Z_i + W_i \cdot Z_i +  Z_i^2$. By the well
known properties of a fundamental cycle, we have $ Z_i^2 = -2$,
and  $W_i \cdot Z_i \leq 0$, therefore $\tilde {\Ga } \cdot Z_i
\geq 2$, and we conclude by the previous inequality that $\tilde
{\Ga } \cdot Z_i  = 2$.

10) By 6) and 8) it follows that if we have two fundamental cycles
which are intersected by $\tilde {\Ga }$,  both are not reduced.
By the standard classification of fundamental cycles, this means
that the corresponding rational double points are not of type
$A_n$, or, equivalently, that on the fibre $F \cong \PP^1$ of
which $\Ga$ is the inverse image, we have two triple points. This
however contradicts 1), and shows that one of the cases [S1] or
[S2] occurs.

11) Let us consider then the former case in 8), where $\tilde {\Ga
} \cdot \tilde {Z_1} = 4$, and there is only one fundamental cycle
which is intersected by $\tilde {\Ga }$, so we have $\tilde {Z_1}
= \tilde {Z} $ and we may write accordingly $Z$ for the
fundamental cycle, and $Z' = Z_{red}.$  Since $ \tilde {Z} ^2 =
-4$, $Z^2= -2$, we can write as in 9) $ \tilde {Z} = Z + W$, and $
-4 = \tilde {Z} ^2 = Z^2 + W^2 + 2 W \cdot Z$, and we get a sum of
non positive terms, where the first two are even and strictly
negative. Hence follows that $ -2 = W^2 , \ W \cdot Z = 0 , \tilde
{\Ga } \cdot Z = \tilde {\Ga } \cdot W = 2$ (note that
$0 = Z \cdot K_{S_0} = Z \cdot (\tilde{\Ga} + Z + W)$).

Thus again the
fundamental cycle corresponds to a triple point of the branch
curve, and $\tilde {\Ga } $ intersects $Z'$ in a smooth point,
belonging to a (-2)-curve $A$ which appears with multiplicity $2$
in both $Z$ and $W$. Write $W = \sum r_i A_i$ with $ r_i \geq 0$,
then $A_i\cdot Z=0$ for all $i$. So the equation $A_i\cdot
(\tilde{\Gamma}+Z+W)=0$ implies $A_i\cdot W=-A_i\cdot
\tilde{\Gamma}\le 0$. Also we have seen that $W$ is a sum of only
those $ A_i$'s which are orthogonal to $Z$. Moreover, since the
point $ P = A \cap \tilde {\Ga }$ is invariant under the
involution $\sigma$, we see that $A$ is pointwise $\sigma$-fixed.
Indeed, both $A$ and $\tilde{\Ga}$
are $\sigma$-stable and their tangents are eigenvectors
of the action $\sigma_*$ on the tangent space at $P$, but
$\tilde{\Ga}$ is not pointwise fixed and if also $A$ were not  we would have
an isolated fixed point, a contradiction.

\par
Thus after we divide by the involution we
obtain a (-4)-curve $E$, image of $A$, such that $\tilde {\Ga }$
is the inverse image of a transversal curve $\tilde {F }$ meeting
$E$ precisely in the point $p$ image of $P$.

12) Let us analyse this last case in terms of  the double covering
$\Ga \ra F$, where $ F \cong \PP^1$. Since, on $S_0$,  $\tilde
{\Ga }$ is smooth of genus $0$, it follows that this covering is
branched on the point $P_{\infty}$ and on another point $P \in F
\cap B$, where  the branch locus $B$ of the double covering meets $F$ with
intersection multiplicity $(B \cdot F)_P = 5$ (observe also that $B$
 does not contain $F$
as a component, else $K_{S_0} \ge 2 \tilde{\Ga}$, absurd.)

Because $Y$ has only Rational Double Points as 
singularities,the branch curve $B$ of the  double covering
 has only simple singularities (see \cite{BPV}). Since the
fundamental cycle is not reduced, then $B$ has a triple point at
$P$. After blowing up $P$ we get a (-1)-curve $E_1$ and the full
transform of $F$ is then $ F' + E_1$, and the new branch locus is
$ B' + E_1$, where $B'$ is the proper transform of $B$. We know
that the curve $E$ occurring in the normal crossing resolution of
the branch locus has multiplicity $2$ in the full transform
of $F$ (since $A$ has multiplicity $4$ in
$D = \tilde{\Ga} + \tilde{Z}$), 
so $E$ cannot be the proper transform of $E_1$, and the new
branch locus has a point of multiplicity $3$ at the point $P' = F'
\cap E_1$, and we must blow up $P'$, obtaining a (-1)-curve $E_2$
which is part of the new branch locus, together with the strict
transform $B''$ of $B'$ and the strict transform $E'_1$ of $E_1$.

Since we had $ B' \cdot F' = 2$, and $ B'$ has multiplicity $2$ at
$P'$, it follows that now $F'' \cap B'' = \emptyset$. Moreover
also $F'' \cap E'_1 = \emptyset$, therefore $E$ is the strict
transform of $E_2$. Since  $E$  has self intersection equal to
$-4$ we need three further blow ups of points (possibly infinitely
near) on $E_2$.

13) Since the proper transform $B'$ is singular we can exclude
that we have a rational double point of type $E_6$ or of type
$E_7$. The other two cases are separated accordingly as follows
(see \cite{BPV}, II (8.1) and III (7.1) for the one-to-one
correspondence between the type of curve singularity 
of the branch locus $ B \subset Q_2$  of the double cover
$Y \rightarrow Q_2$ and the type of surface singularity
 at the corresponding point on the canonical model $Y$).

{\bf [S3.1]} $ (B' \cdot E_1)_{P'} = 2$ implies that we have the $D_n$
case, since $B$ has then two distinct tangents at $P$.

{\bf [S3.2]} $ (B' \cdot E_1)_{P'} = 3$ implies that we have the $E_8$
case, since $B$ has then only one tangent at $P$.

14) In both cases, we observe that $\tilde{Z}$ is the pull-back of
$E'_1  + 2 E_2$, i.e., the pull back of the maximal ideal of $P$
plus the pull back of the maximal ideal of $P'$. Moreover, in case
[S3.2], there is only one (-2)-curve $A$  which occurs in $Z$ with
multiplicity two, and such that $ A \cdot Z = 0$. In case [S3.1] we
see instead that $A$ is the  curve corresponding to the vertex at
distance three from the asymmetrical end (observe that our
assumptions imply $ n \geq 6$).

15) We proceed by observing that it suffices to verify the
statement for one blow up of $S_0$ where we have normal crossings
for $ C : = f^* (D)$.

\begin{lem}\label{nc}
Let $C \subset S$ be a normal crossing divisor, and let $ g : S'
\ra S$ the blow up of a point $P$. Then,  if we set $ C' = g^*
(C)$, and $ \Theta_t : =\roundup {t C}$, $ \Theta_t' : =\roundup {
t C'}$, then  $ h^0 ( K_S + \Theta_t) =  h^0 ( K_{S'} + \Theta'
_t)$, , for all $t \in (0, 1)$.

More generally, let $ C = \sum_i n_i C_i$ be the decomposition of
$C$ as a sum of irreducible analytic branches at $P$, and let $m_i
: = mult_P (C_i)$, then the above equality holds if there are two
smooth local branches, or just one branch ( i.e., $n_1=1, n_j = 0
\ \forall j \geq 2$) of multiplicity $m=2$, provided $ t \in
(0,1)$.
\end{lem}

\begin{proof} Let $E$ be the exceptional divisor of $g$, let
$ C = \sum_i n_i C_i$ be the decomposition of $C$ as
a sum of irreducible divisors: then
$$  C' = g^* (C)   = \sum_i n_i C'_i +  \sum_i n_i m_i E,$$
where $ C'_i$ is the proper transform of $C_i$, and $m_i : =
mult_P (C_i)$.

Taking the round up, we obtain
\begin{align*}
  \Theta' _t  &=\roundup { t  C'} =  \sum_i \roundup { t n_i} C'_i
+ \roundup { t \sum_i m_i n_i}E \\
&= g^* (\Theta_t)  +  (  \roundup { t \sum_i m_i n_i} - \sum_i
\roundup { t n_i} m_i )  E .
\end{align*}

Since $ K_{S'} = g^* (K_S)+E$, $ K_{S'} + \Theta' _t = g^* (  K_S
+ \Theta_t)  + [ 1 + \roundup { t \sum_i m_i n_i} - \sum_i
\roundup { t n_i} m_i ]  E$ and it suffices that the integer in
the square brackets is non negative in order to conclude the
desired equality. Notice that the calculation is entirely local,
so that we can replace the global decomposition by the local
decomposition in analytic branches.

The normal crossings case is a special case of the one where all
the multiplicities satisfy $ m_i =1$ : in this case we want the
inequality
\begin{align*}
\tag{**} 1 + \roundup { \sum_i  t n_i} \geq \sum_i \roundup { t
n_i}\end{align*}
to hold. This is obvious if there are exactly two
terms, since for any two real numbers $a,b$ holds $ 1 + \roundup{
a +b} \geq \roundup{ a }+ \roundup{ b}$.

For only one branch, we want $ 1 + \roundup { tm} \geq \roundup
{t} m$, and this is true for $m=2$, since $t <1$.

\end{proof}

{\bf Case [S1]}: we take $C = \tilde {\Ga }$: it is irreducible of
arithmetic genus equal to $p \in \{ 1, 2 \}$, therefore, for each
$t \in (0,1)$ $C$ is equal to the round up of $tC$, and $ h^0 (\om_C )
= p$. If  $C$ has normal crossings, we are done by the previous
lemma (choose $D' = C$ in  2)).

Assume the contrary and assume first  $ p=1$: then $C$ has an
ordinary cusp, thus the hypothesis of the lemma above applies.
After a blow up we get two  smooth tangent  branches ( $ n_1=1,
n_2 =2$), and the lemma still applies. We then get  three smooth
transversal branches where $ n_1=1, n_2 =2, n_3 = 3$: the
inequality (**) holds, provided $ 1/6 < t \leq 1/3$ (since it is
equivalent then to  $\roundup {6 t} \geq 2$) and we are done,
since after this blow up we get global normal crossings for the
full transform.

Assume now that $C$ does not have normal crossings, and that $ p=2$: we
have just verified that an ordinary cusp gives no problem ( as
well as a node). We use now the fact that $C$ has only double
points as singularities, so we have to verify that a tacnode $ y
^2 = x^4$ and a higher cusp $ y ^2 = x^5$ give no problem (higher
singularities are excluded by  point 1)).

For a tacnode we get two smooth branches, so the lemma applies,
and after the first blow up we get three smooth transversal
branches , with $ n_1=1, n_2 =1, n_3 = 2$, thus (**) applies again
if $ 1/4 < t \leq 1/3$ (since (**) is equivalent then to $\roundup
{4 t} \geq 2$)  and after this blow up we get normal crossings.

In the case of the higher cusp, we get one branch of multiplicity
$2$, so the lemma applies; after the first blow up we get a
reduced  ordinary cusp transversal to a smooth branch, occurring
with multiplicity $2$. In this case we have to verify that $ 1 +
\roundup {4 t} \geq  \roundup {2 t} + 2 \  \roundup { t}$ , but
this clearly holds for  $ 1/4 < t \leq 1/3$.

After a further blow up, we get two smooth branches, tangent, and
with $ n_1=1, n_2 =4$, so the lemma applies. A further blow up,
the last before we get normal crossings, yields a point where
three smooth branches meet transversally, and $ n_1=1, n_2 =4, n_3
= 5$: we have to verify whether (**) holds, i.e., $  1 + \roundup
{ 10 t} \geq  \roundup { t} + \roundup{ 4 t} + \roundup {5t}$. But
this holds clearly for $ 3/10 <  t \leq 1/3$ (else , for $1/5 <  t
\leq 3/10$ there is a loss by $1$, which however would not trouble
us since we started with $ p=2$, and we only want the arithmetic
genus above to be at least $1$).
We now treat the remaining cases one by one,  using  2).

{\bf Case [S2.1]}: $D'$ already has normal crossings, and is
reduced, thus there is nothing to prove.

{\bf Case [S2.2]}: here $D'$ consists of two smooth tangent
divisors $ \cong \PP^1$, so its arithmetic genus is $ p=1$. This
is exactly the case of the tacnode, which we already treated, thus
this case is also settled.

{\bf Case [S2.3]}: $D$ ($\sim K_{S_0}$)  now has arithmetic genus
2, and does not have normal crossings exactly at the point where $A,
A', \tilde{\Ga} \cong \PP^1$ meet transversally. The local
multiplicities are $ 2,2,1$, thus for $ 1/5 < t \leq 1/3$ we
obtain $  1 + \roundup { 5 t} \geq  \roundup { t} + \roundup{ 2 t}
+ \roundup {2t}$. Thus we are done as in the above lemma.

{\bf Case [S2.4]}: $D'$ already has normal crossings  and has
arithmetic genus 1. So there is nothing to prove.

{\bf Case [S3.1]}: an explicit calculation, probably well known,
(cf. \cite{BPV}, page 65, lines 2-3) shows that the full transform
of the maximal ideal of $P$ is the fundamental cycle $Z$ of $D_n$,
while the full transform of the maximal ideal of $P'$, which is
then $W$, is the fundamental cycle of the $D_{n-2}$ configuration
obtained by deleting the asymmetric end and its neighbour.

In this case let us choose as $D' = \tilde{\Ga} + 2 W < D =
\tilde{\Gamma} + \tilde{Z}$:
since all the multiplicities of the components of $D'$ are then 
 either $1$, $ 2$ or $ 4$, it follows that
for $ 1/4 < t \leq 1/3$ the round up $\Delta_t' : = \roundup{ t
D'}$ equals $\tilde{\Ga} + W$. Since $W^2 = -2$, $\tilde{\Ga}
\cdot W = 2$, the self intersection $(\Delta_t') ^2 = -1$, thus
$\Delta_t' $ has arithmetic genus $1$
(topologically it is of elliptic type $D_{n-2}^*$ or $I_{n-6}^*$
in Kodaira's notation) and this case is settled
by virtue of 2).

{\bf Case [S3.2]}:  Also in this case an explicit calculation,
probably well known, (cf. \cite{BPV}, page 65, lines 2-3) shows
that the full transform of the maximal ideal of $P$ is the
fundamental cycle $Z$ of $E_8$, while the full transform of the
maximal ideal of $P'$, which is then $W$, is the fundamental cycle
of the $E_7$ configuration obtained by deleting the furthest end.

In this case we write the multiplicities for the components of
$\tilde{Z}$ starting from left to right (i.e., from middle length
end (i.e., $A$) to longest end), and then we give the multiplicity
for the shortest end: we get the sequence $4,7,10,8,6,4,2$ and then $5$.
We may choose for convenience $D'$ as $\tilde{\Ga} +\tilde{Z}$
   minus twice the
longest end and minus its neighbor, i.e., we change the sequence
to $4,7,10,8,6,3,0, 5$. If we now choose $ 3/10 < t \leq 1/3$, one
can easily calculate that the round up $\Delta_t' : = \roundup{ t
D'}$ equals $\tilde{\Ga} + W$
(topologically, it is of elliptic type $E_7^*$ or $III^*$
in Kodaira's notation), and we are done as in the previous
case.
\end{proof}

\section{\bf The Noether inequality}

\begin{thm}\label{inequality} Let $X$ be a minimal Gorenstein 3-fold of
general type with canonical singularities. Assume either
$p_g(X)\le 2$ or that $|K_X|$ is composed with a pencil of
surfaces of type (1,2). Then
$$K_X^3\ge \frac{7}{5}p_g(X)-2.$$
\end{thm}

\begin{proof} As we have seen in \ref{notation}, we may take $X$ to be
factorial with only terminal singularities. Because $K_X^3\ge 2$, the
inequality is automatically true for $p_g(X)\le 2$.

We may suppose, from now on, that $p_g(X)\ge 3$. Denote by
$f:X'\longrightarrow B$ the  fibration induced from $\Phi_{|K_X|}$.
Let $S$ be a general fiber of $f$. Lemma 4.5 of
\cite{chenjmsj2004} says $0\le b=g(B)\le 1$. Theorem \ref{cat}
says that we may take $\alpha_3=\frac{3}{10}$ for a general fiber
$S$ of $f$; see the first part of 2.3.
\smallskip

{\bf Case 1}. $b=1$. We may write $a=p_g(X)=10m+c$ where $m \geq 0$ and
$0\le c\le 9$ (obviously, for $m=0$ we have $3\le c\le 9$).

When $m>0$, we take $m_3 : =7m+\Box_c$, where $\Box_c : =-1$, 0, 1, 2, 2,
3, 4, 4, 5, 6 respectively when $0\le c\le 9$. Then one sees
that
$$1-\frac{m_3}{a}>\alpha_3=\frac{3}{10}.$$
Take $m_2=a-1-m_3=3m+\triangledown_c$, where
$\triangledown_c : = c-1-\Box_c$. Then the inequality (2.5) gives
\begin{align*}
K_X^3&\ge 6(7m+\Box_c)+4(3m+\triangledown_c)-4p_g(X)+4\\
&=54m-4p_g(X)+6\Box_c+4\triangledown_c+4\\
&=\frac{7}{5}p_g(X)-\frac{7}{5}c+2\Box_c\\
&\ge \frac{7}{5}p_g(X)-2.
\end{align*}
When $m=0$, we have $3\le a=c\le 9$. Take $m_3=$2, 2, 3, 4, 4, 5,
6 respectively when $3\le c\le 9$. We may easily check that
$1-\frac{m_3}{a}>\frac{3}{10}$. Take $m_2 : = a-1-m_3=$0, 1, 1, 1, 2,
2, 2 respectively when $3\le c\le 9$. By inequality (2.5), we may
verify case by case that $K_X^3
>\frac{7}{5}p_g(X)-2. $

{\bf Case 2}. $b=0$. We may write $a=p_g(X)-1=10m+c$ where $m \geq 0$
and $0\le c\le 9$ (for $m=0$ we have $2\le c\le 9$).

Again when $m>0$, we take $m_3 : =7m+\Box_c$
where $\Box_c=-1$, 0, 1, 2, 2, 3, 4, 4, 5, 6 respectively when
$0\le c\le 9$. Then the calculation is similar to Case 1. Take
$m_2=a-1-m_3=3m+\triangledown_c$ where
$\triangledown_c=c-1-\Box_c$. Then the inequality (2.5) gives
\begin{align*}
K_X^3&\ge 6(7m+\Box_c)+4(3m+\triangledown_c)-4p_g(X)+
4h^2({\mathcal O}_X)+10\\
&\ge 54m-4p_g(X)+6\Box_c+4\triangledown_c+10\\
&=\frac{7}{5}p_g(X)-\frac{7}{5}c+2\Box_c+\frac{3}{5}\\
&\ge \frac{7}{5}p_g(X)-\frac{7}{5}.
\end{align*}
When $m=0$ and $2\le a=c\le 9$, one can  in a similar way 
verify that $K_X^3>\frac{7}{5}p_g(X)-\frac{7}{5}. $
\end{proof}

\begin{setup}{\bf Proof of the main results.}\end{setup}
\begin{proof}

Now both Theorem \ref{main} and Corollary \ref{application} follow
directly from \ref{known} and Theorem \ref{inequality}.
\end{proof}

\begin{rem} A quite natural problem left to us is the possibility
of  generalizing  Theorem \ref{main} to the case where $X$
is Gorenstein minimal. Unfortunately the method of Theorem 4.3 of
\cite{chenjmsj2004} only works when $X$ is smooth. One needs a new
method to treat the  difficult case where $X$ is canonically
fibred by curves of genus 2. However we would like to put forward
the following:
\end{rem}

\begin{conj}  The Noether inequality 
$$K^3\ge \frac{2}{3}(2p_g-5)$$
holds for any projective minimal Gorenstein 3-fold of
general type $X$.
\end{conj}
\bigskip

\noindent{\bf Acknowledgment.} This paper was begun while Zhang
was visiting Fudan University at the end of 2004. It was finished
when Chen was visiting both Universit$\ddot{\text{a}}$t Bayreuth
and Universit$\ddot{\text{a}}$t Duisburg-Essen in the summer of
2005 financially supported by the joint Chinese-German project on
"Komplexe Geometrie" (DFG \& NSFC). The authors would like to
thank these 3 universities for their support. Especially Chen
would like to thank Ingrid Bauer, H\'el\`ene Esnault, Thomas
Peternell, Eckart Viehweg for their hospitality and effective
discussion.


\end{document}